\documentclass{amsart}
\usepackage{mathrsfs}
\usepackage{graphics}
\usepackage{color}
\usepackage{isolatin1}
\title[Rational semimodules over the max-plus semiring]{Rational semimodules over the max-plus semiring and geometric approach of discrete event systems}
\author{St\'ephane Gaubert}
\address{INRIA, Domaine de Voluceau, 78153, Le Chesnay C\'edex, France}
\email{Stephane.Gaubert@inria.fr}
\author{Ricardo Katz}
\address{CONICET. Postal address: Dep. of Mathematics, Universidad Nacional de Rosario, Avenida
Pellegrini 250, 2000 Rosario, Argentina}
\email{rkatz@fceia.unr.edu.ar}

\date{November 8, 2002}
\newtheorem{theorem}{Theorem}[section]
\newtheorem{proposition}[theorem]{Proposition} 
\newtheorem{lemma}[theorem]{Lemma}
\theoremstyle{definition}
\newtheorem{definition}[theorem]{Definition} 
\theoremstyle{remark}
\newtheorem{remark}[theorem]{Remark}
\newtheorem{example}[theorem]{Example}
\newcommand{\NEW}[1]{{\em #1}}
\newcommand{\zero}{\varepsilon}
\newcommand{\unit}{e}
\newcommand{\cG}{\mathcal{G}}
\newcommand{\cB}{\mathcal{B}}
\newcommand{\cO}{\mathcal{O}}
\newcommand{\cP}{\mathcal{P}}
\newcommand{\cR}{\mathcal{R}}
\newcommand{\sR}{\mathscr{R}}
\newcommand{\cS}{\mathcal{S}}
\newcommand{\cW}{\mathcal{W}}
\newcommand{\cX}{\mathcal{X}}
\newcommand{\cY}{\mathcal{Y}}
\newcommand{\cZ}{\mathcal{Z}}
\newcommand{\set}[2]{\{#1\mid\,#2\}}
\newcommand{\Bigset}[2]{\big\{#1\mid\,#2\big\}}
\newcommand{\mrm}[1]{\text{\rm #1}}

\newcommand{\diag}{\mrm{diag}}
\newcommand{\im}{\mrm{Im}\,}
\renewcommand{\ker}{\mrm{Ker}\,}
\newcommand{\Ext}{\mrm{Ext}}
\newcommand{\Hom}{\mrm{Hom}}
\newcommand{\R}{\mathbb{R}}
\newcommand{\N}{\mathbb{N}}
\newcommand{\Z}{\mathbb{Z}}
\newcommand{\Q}{\mathbb{Q}}
\newcommand{\Rmat}{\mathsf{R}}
\newcommand{\Omat}{\mathsf{O}}
\newcommand{\Zpm}{\Z\cup\{\pm\infty\}}
\newcommand{\Zm}{\Z\cup\{-\infty\}}
\newcommand{\Rm}{\R\cup\{-\infty\}}
\newcommand{\Zp}{\Z\cup\{+\infty\}}

\newcommand{\pattern}{\{0,\pm\infty\}}
\newcommand{\rmax}{\R_{\max}}

\newcommand{\qmax}{\Q_{\max}}
\newcommand{\nmin}{\N_{\min}}
\newcommand{\nmax}{\N_{\max}}
\newcommand{\nmaxb}{\bar{\N}_{\max}}
\newcommand{\zmin}{\Z_{\min}}
\newcommand{\zmax}{\Z_{\max}}
\newcommand{\zmaxb}{\bar{\Z}_{\max}}
\newcommand{\vect}{\mrm{span}\,}
\newcommand{\smallvect}{\mrm{span}}
\newcommand{\mynat}[1]{\mrm{Nat}(#1)}
\begin{document}
\begin{abstract}
We introduce rational semimodules 
over semirings whose addition is idempotent,
like the max-plus semiring,
in order to extend the geometric approach
of linear control to discrete event systems.
We say that a subsemimodule
of the free semimodule $\cS^n$ over a semiring $\cS$
is rational if it has a generating family
that is a rational subset of $\cS^n$,
$\cS^n$ being thought of as a monoid under
the entrywise product. We show
that for various semirings of 
max-plus type whose elements are integers,
rational semimodules
are stable under the natural algebraic
operations (union, product, direct and inverse image, intersection, 
projection, etc). 
We show that the reachable
and observable spaces of max-plus linear
dynamical systems are rational, and give
various examples. 
\end{abstract}
\maketitle
\tableofcontents
\section{Introduction}
In this paper, we define a new class of
semimodules over max-plus type semirings,
that we call \NEW{rational semimodules},
and study their properties.

This work is motivated by the max-plus algebraic
approach of discrete event systems.
It is well known (see in particular~\cite{cohen89a,bcoq,idem,ccggq99})
that a subclass of discrete event systems subject
to synchronization constraints, comprising
examples of manufacturing systems, transportation
networks, and computer networks, can be
modeled by max-plus linear dynamical systems.
An open question (see~\cite{ccggq99})
is to develop the analogue of Wonham's
geometric approach~\cite{wonham} for the control
of max-plus linear dynamical systems.
As in classical linear system theory,
many control problems can be
phrased in terms of semimodules
(semimodules over semirings
are defined like modules over
rings, mutatis mutandis).
A difficulty of this approach, however, 
is that max-plus semimodules have very
different properties from vector spaces.
In particular, a subsemimodule
of a free finitely generated semimodule
need not be free or finitely generated,
and so even the simplest spaces in control theory,
the reachability space and the observability ``space'' or congruence,
need not be finitely generated (see the examples in~\S\ref{subsec-display}).
Therefore, new algebraic tools are needed to ``replace'' the theory
of rank which is so useful in classical linear control,
and effective methods must be designed to handle semimodules
with an infinite number of generators.

Several results are known on max-plus semimodules,
including notions of basis and extremal
points~\cite{moller88,wagneur91,maxplus97,gaubert98n}, 
direct sums~\cite{CGQ96a}, projective semimodules~\cite{CGQ97a},
separation theorems~\cite{zimmermann77,shpiz,litvinov00,gaubert01a}.
However, the issue of computing effectively with 
non finitely generated semimodules
does not seem to have been
raised previously in the literature.

In this paper, we extend the notion of finitely
generated semimodule as follows:
we say that a semimodule $\cX\subset \cS^n$
is \NEW{rational}
if it has a set of generators that is a rational
subset of $\cS^n$, where $\cS^n$ is
thought of as a monoid under the entrywise product,
see Definition~\ref{def-rat} below.
Rational sets over monoids,
and in particular,
rational sets of $(\N^k,+)$
or $(\Z^k,+)$, or {\em semilinear} sets,
are well known objects in computer science,
see~\cite{GinSpa,eilenberg69}. The typical
semiring to which our notions apply is
the semiring of max-plus integers,
$\Z\cup\{-\infty\}$, equipped with max
as addition, and the usual addition as 
multiplication: then, up to
technical details related to the infinite
element, rational semimodules are 
semimodules generated by semilinear
sets of $\Z^n$. 

We show that rational semimodules
are closed under the natural algebraic
operations, like union, direct sum, direct
and inverse image, intersection, projection,
and by taking the orthogonal.
Whereas the closure under union, direct sum,
and direct image, can be proved in a natural way, 
our proof of the other properties
relies on Presburger arithmetics,
which leads to expensive algorithms~\cite{oppen}. 
Finding direct, computationally
more efficient proofs, leads
to interesting combinatorial problems.
In fact, even for finitely
generated semimodules, 
algorithmic issues remain
difficult, see Remark~\ref{rk-algo}
below.

The paper is organized as follows. In 
\S\ref{sec-pres},
we recall classical definitions and facts about
rational sets and idempotent semirings, and establish
preliminary results. We extend the definition
of the Presburger logic to a naturally ordered
idempotent semiring $\cS$,
and show, as a slight extension of the theorem of Ginzburg
and Spanier~\cite{GinSpa}, 
that idempotent semirings like $(\Zm,\max,+)$
have the property that the subsets of $\cS^n$ defined by formulas
of the first order logic of $(\cS,e,\otimes, \preceq)$,
where $e$ is the unit, $\otimes$ the product, and $\preceq$
the natural order, are exactly the rational subsets of $\cS^n$.
In~\S\ref{sec-closure}, we use these results to show that
the class of rational semimodules is closed
under various algebraic operations.
In~\S\ref{secAyE}, we illustrate the results
by discrete event systems problems,
and give various examples and counter examples.
We show in~\S\ref{subsec-obs} that max-plus reachable
spaces and observable congruences are rational,
and then, in~\S\ref{subsec-display}, we give examples
of reachable semimodules and observable
congruences. In~\S\ref{subsec-machine}, we illustrate
the reachability and observability
notions by discussing a simple example of manufacturing
systems (three machines in tandem). In~\S\ref{subsec-cex},
we give counter examples showing that the integrity
assumptions that we used are useful, and in~\S\ref{subsec-nc},
we give a counter example showing
that the noncommutative analogues of reachable
spaces need not be rational.

\section{Presburger Logic over Idempotent Semirings}
\label{sec-pres}

Let us recall some definitions and results. Let $(M,\cdot )$ be a monoid, i.e. a set with an associative multiplication
and a two sided unit $1_M$. The class 
of {\em rational}\/ subsets of $M$ is the
least class $\sR$ of subsets of $M$ satisfying the following conditions:
\begin{enumerate}
\item If $U$ is a finite set then $U\in\sR$;
\item If $U$, $V\in\sR$ then $U\cup V\in\sR$;
\item If $U$, $V\in\sR$ then
$U\cdot V=\set{m}{m=u\cdot v,u\in U,v\in V}\in\sR$;
\item If $U\in\sR$ then $U^*=U^0\cup U\cup U^2\cup\cdots\in\sR$,
\end{enumerate}
where by convention $U^0=\{1_M\}$.
A subset $U\subset M$ is called \NEW{semilinear} if it can be written as a finite union of sets of the form $\{x\} \cdot B^{*}$,
where $x\in M$ and $B$ is a finite subset of $M$. We 
shall use throughout the paper the following
classical result (see~\cite{eilenberg69}):
in a commutative monoid, rational and
semilinear subsets coincide. 
 
A \NEW{semiring} is a set ${\cS}$ equipped with two internal composition laws $\oplus $ and $\otimes $, called addition and multiplication respectively,
such that ${\cS}$ is a commutative monoid for addition,
${\cS}$ is a monoid for multiplication,
multiplication distributes over addition,
and the neutral element for addition is absorbing for multiplication.
We will denote by $\zero $ the neutral element for addition and by $e$ the neutral element for multiplication. We will
sometimes denote by $(\cS, \oplus,\otimes)$
or $(\cS,\oplus,\otimes,\zero,e)$ the semiring $\cS$.
In this paper, we are mostly interested in the max-plus semiring $\rmax$, 
which is the set $\R\cup\{-\infty\}$ equipped with
$\oplus = \max $ and $\otimes =+$. 
The semiring $\rmax$ is \NEW{idempotent}: 
$x\oplus x= x$ for all $x\in\cS$.
An idempotent semiring $(\cS,\oplus,\otimes)$ is equipped
with the \NEW{natural order}
$\preceq$, which is defined by:
$$x\preceq y \Leftrightarrow x\oplus y=y.$$
With this order, $x\oplus y$ is the least upper bound of the set
$\{x,y\}$ (see \cite{bcoq}). 

Several variants of the max-plus semiring $\rmax$ can
be found in the literature. Indeed, to any submonoid $(M,+)$ 
of $(\R,+)$ is associated a semiring with set
of elements $M\cup\{-\infty\}$, and laws $\oplus =\max, \otimes =+$.
We denote this semiring by $M_{\max}$. Symmetrically,
the semiring $M_{\min}$ is the set $M\cup\{+\infty\}$,
equipped with $\oplus=\min$ and $\otimes =+$.
For instance, taking $M=\N$, we get $\nmin=(\N\cup\{+\infty\},\min,+)$, a semiring known as the \NEW{tropical semiring}
after the work of Simon~\cite{simon78}
(see~\cite{pin95} for a recent overview). The semiring
$\zmin=(\Zp,\min,+)$ is sometimes
called the {\em equatorial semiring}~\cite{krob93b,krob94a}.
One can also add a maximal element (for the natural order) 
to the semirings
$M_{\max}$ and $M_{\min}$: this yields
the semirings ${\bar{M}}_{\max}=(M\cup\{\pm\infty\},\max,+)$
and ${\bar{M}}_{\min}=(M\cup\{\pm\infty\},\min,+)$.
Since the zero element is $\zero=-\infty$
in $\bar{M}_{\max}$
and $\zero=+\infty$
in $\bar M_{\min}$, in these semirings,
the value of $(-\infty)+(+\infty)=(+\infty)+(-\infty)$
is determined by the rule
$\zero\otimes x=x\otimes \zero=\zero$.

It is convenient to formalize the class of semirings
to which our results apply, 
by extending the classical definition of Presburger logic,
as follows. We refer the reader
to~\cite{GinSpa,bes} for more information
about Presburger logic. Our presentation
follows~\cite{GinSpa}.

Let $(\cS,\oplus,\otimes,\zero,e)$
be an idempotent commutative semiring with natural order $\preceq$. 
We consider formulas or statements about the elements of $\cS$. 
The set $\cP$ of \NEW{first-order logic formulas} of
$(\cS,e,\otimes ,\preceq)$
is by definition
the smallest class of formulas satisfying
the following five conditions: 
\begin{enumerate}
\item For any nonnegative integers $k_i,r_i,1\leq i \leq n$, 
\begin{align}
\label{p1}
\bigotimes _{i=1}^{n} x_i^{k_i} \preceq \bigotimes _{j=1}^{n} x_j^{r_j}
\end{align}
is a formula in $\cP$. Here $x_i^{k_i}$ denotes 
$x_{i}\otimes \cdots \otimes x_{i}$, where $x_i$ is repeated $k_i$ times,
and we adopt the convention $x_i^0=e$.
The free variables of this formula are $x_1,\ldots ,x_n$;
\item If $P_1$, $P_2$ are in ${\cP}$, so is their conjunction $P_1\wedge P_2$. The set of free variables of $P_1\wedge P_2$ is the set of free variables of $P_1$ union the set of free variables of $P_2$;
\item If $P_1$, $P_2$ are in ${\cP}$, so is their disjunction $P_1\vee  P_2$. The set of free variables of $P_1\vee P_2$ is the set of free variables of $P_1$ union the set of free variables of $P_2$;
\item If $P$ is in ${\cP}$, so is its negation $\neg P$. The free variables of $\neg P$ are the free variables of $P$.
\item If $P( x_1,\ldots ,x_n)$ is in $\cP$ 
and has the free variables $x_1,\ldots,x_n$,
then for each $1\leq i\leq n$, 
the formula $( \exists x_i) P(x_1,\ldots,x_n)$
is in $\cP$ and its free variables are $x_j$ for
$1\leq j \leq n$ and $j\neq i$.   
\end{enumerate}
In the sequel, we will simply call a formula of $\cP$
a \NEW{Presburger formula} of $(\cS,e,\otimes,\preceq)$.
\begin{remark}\sloppy
If $P( x_1,\ldots ,x_n)$ is in ${\cP}$, then for each $1\leq i \leq n$, the formula $( \forall x_i) P( x_1,\ldots ,x_n) $ 
is regarded as a first-order logic formula of
$( {\cS}, e,\otimes ,\preceq) $ because it is equivalent to 
$$\neg \left( \exists x_i \right) \left( \neg P\left( x_1,\ldots ,x_n \right) \right) .$$
Similarly if $P$ and $Q$ are in ${\cP}$, then $P\Rightarrow Q$ is regarded as a first-order logic formula of $\left( {\cS}, e,\otimes ,\preceq \right) $.
\end{remark}
\begin{remark}
The formula $x_{n+1}=\bigoplus_{i=1}^n x_i$ is regarded as a first-order logic formula of $\left( {\cS}, e,\otimes ,\preceq \right) $ because it is equivalent to 
\begin{align*}
&\left( x_1\preceq x_{n+1} \right) \wedge \ldots \wedge \left( x_n\preceq x_{n+1} \right) \wedge \\
& \qquad \left\{ ( \forall x_{n+2}) \left[ \left( \left( x_1\preceq x_{n+2} \right) \wedge \ldots \wedge \left( x_n\preceq x_{n+2} \right) \right) 
\Rightarrow x_{n+1}\preceq x_{n+2} \right] \right\} .
\end{align*}
\end{remark}
For readability, we will allow the use of arbitrary letters
(rather than $x_1,x_2,\ldots$) for the variables of formulas,
so that we will regard for instance $y=\bigoplus_{i=1}^n x_i$
as a Presburger formula with free variables $x_1,\ldots,x_n, y$.

We say that a subset 
$D\subset\cS^n$ is \NEW{definable} in the first-order logic of
$(\cS,e,\otimes,\preceq)$ if there exists a formula
$P(x_1,\ldots ,x_n)$ in $\cP$, with $n$ free variables 
$x_1,\ldots ,x_n$, such that: 
$$D=\set{( x_1,\ldots ,x_n)\in\cS^n}{P( x_1,\ldots ,x_n) \mrm{ is true } }\enspace .$$

\begin{definition}
An idempotent commutative semiring $\left( {\cS},\oplus ,\otimes , e,\zero \right) $ has the \NEW{Presburger property} if the subsets of ${\cS}^n$ definable in the first-order logic of $\left( {\cS}, e,\otimes ,\preceq \right) $ are precisely the rational sets of $\left( {\cS}^n,\otimes \right)$.
\end{definition}
We shall need the following extension
of the theorem of Ginzburg and Spanier (see~\cite{GinSpa}),
which states that the rational subsets of $(\N^n,+)$
are precisely the subsets definable in the classical
Presburger arithmetics.
\begin{theorem}\label{th-pres}
The idempotent semirings 
$\zmax=(\Zm,\max ,+)$, 
$\zmaxb=(\Zpm,\max ,+)$,
$\nmax=(\N\cup\{ -\infty\} ,\max ,+)$,
$\nmaxb=(\N\cup \{ \pm\infty\} ,\max ,+)$,
and $\nmin=(\N\cup\{+\infty\} ,\min ,+)$ 
all have the Presburger property.
\end{theorem}
The proof of Theorem~\ref{th-pres}
relies on the following:
\begin{lemma}\label{lem-new}
A subset $R\subset ((\Zpm)^n,+)$ is
rational if and only if it can be written
as a finite union of sets of the form
\begin{align}
\label{pfinite}
\{a\}+ \{\bar r_1,\ldots,\bar r_k\}^* \enspace,
\end{align}
where $a\in (\Zpm)^n$ and $\bar r^1,\ldots,\bar r^k\in \Z^n$.
\end{lemma}
\begin{proof}
Using the characterization of rational sets
as semilinear sets, it suffices to show
that any set $\{b\}+\{r^1,\ldots,r^k\}^*$,
with $b,r^1,\ldots,r^k\in(\Zpm)^n$,
can be rewritten as a finite union of sets
of the form~\eqref{pfinite}. 
Recall the following classical {\em rational identities}
\begin{align}
\label{erat}
(A\cup B)^*&=A^*+B^* \\
\label{erat2}
A^* &= A^0\cup (A+A^*) 
\end{align}
(identity~\eqref{erat} holds for all subsets
$A,B$ of a commutative monoid,
whereas~\eqref{erat2} holds for subsets
$A$ of arbitrary monoids, see e.g.~\cite{conway71}
for more details on rational identities). 
Using~\eqref{erat} and~\eqref{erat2}, we can write:
\begin{align}
\nonumber
\{b\}+\{r^1,\ldots,r^k\}^{*}&=
\{b\}\cup\bigcup_{1\leq i\leq k} \{b+r^i\}+\{r^1,\ldots,r^k\}^{*}\\
&=\{b\}\cup
\bigcup_{1\leq i\leq k}
\left(\{b+r^i\}+\{r^i\}^*+\set{r^j}{1\leq j\leq k,\;j\neq i}^*\right)
\enspace .
\label{e-useful}
\end{align}
Now, for all $r\in (\Zpm)^n$, 
denote by $\bar r\in \Z^n$ the vector
obtained by replacing
infinite coordinates of $r$ by an arbitrary
finite value (say $0$). One
easily gets, using the fact that $(-\infty)+ x=-\infty$
for all $x\in \zmaxb$, and $(+\infty) +x=+\infty$,
for all $x\in \Z$, that 
\begin{align}
\{b+ r\}+\{r\}^*=\{b+ r\}+\{\bar r\}^*\label{e-subs} \enspace .
\end{align}
Using several times the transformations \eqref{erat},~\eqref{e-useful},~\eqref{e-subs},
we express $\{b\}+\{r^1,\ldots,r^k\}^*$
as a finite union of sets
of the form $\{a\}+\{\bar r^1,\ldots,\bar r^k\}^*$,
with $a\in (\Zpm)^n$.
\end{proof}
\begin{proof}[Proof of Theorem~\ref{th-pres}]
By comparison with the result of Ginzburg
and Spanier, the only new difficulty is to
take care of the $\pm\infty$ elements.
Let us consider the case of $\zmaxb$.
The other cases can be proved in the same way.

As a preliminary result, we first check that
every element of $\Zpm$ 
is definable by a Presburger formula of
$(\Zpm,0,+,\leq)$. 
We can regard $y=+\infty$ and $y=-\infty$
as Presburger formulas,
since $y=+\infty$ (resp. $y=-\infty$)
is equivalent to $(\forall x) ( x\leq y)$
(resp. $( \forall x) ( x\geq y)$).
Similarly, $y=1$, which is equivalent
to $\neg( y\leq 0) \wedge(\forall x)(\neg(x\leq 0)\Rightarrow y\leq x)$, will be seen as a Presburger formula.
We note that in $\zmaxb$, the inequality~\eqref{p1}
becomes:
\begin{align}
\label{p11}
\sum _{i=1}^{n} k_i x_i \leq \sum _{j=1}^{n} r_j x_j
\enspace .
\end{align}
Therefore, for any positive integer $r$, 
$y=r$, which is equivalent to $(\exists z)(z=1)\wedge (y=rz)$,
will be seen as a Presburger formula,
as well as $y=-r$, which 
is equivalent to $(\exists z)(z=r)\wedge (0=y+ z)$.
Finally, we denote by $\mynat y$
the Presburger formula $(y\geq 0)\wedge \neg (y\geq +\infty)$,
which expresses the property that $y$ is a natural number.

We next show that every semilinear set of 
$((\Zpm)^n,+)$ is definable 
by a Presburger formula of $(\Zpm,0,+,\leq)$.
Since the family of sets definable
in the first-order logic of $(\Zpm,0,+,\leq)$
is closed under union,
thanks to Lemma~\ref{lem-new},
it is enough to show the following:

(Claim C):
{\em  For all $a\in(\Zpm)^n$, and $\bar r^1,\ldots,\bar r^k\in \Z^n$,
the set~\eqref{pfinite} is definable by a 
formula of the first-order logic of $(\Zpm,0,+,\leq)$.}

Indeed, for each $1\leq i\leq n$ let us define the set
$J_i=\set{1\leq j \leq k}{r^j_i<0}$.
Then, the vector $(x_1,\ldots ,x_n)\in (\Zpm)^n$ belongs
to $\{a\}+\{\bar r^1,\ldots ,\bar r^k\} ^{*}$ if and only if
\[
(\exists y_1),\ldots ,( \exists y_k)
\big(
\mynat{y_1}\wedge \cdots\wedge \mynat{y_k}\wedge \bigwedge\limits_{1\leq i\leq n} P_i(x_i,y_1,\ldots ,y_k)
\big)\enspace,
\]
where: 
\begin{align}
P_i(x_i,y_1,\ldots ,y_k) &= ( \exists z_i) \Big( (z_i=a_i)\wedge
\big( x_i+\sum_{j\in J_i}( -\bar r^j_i) y_j=z_i+\sum_{j\not\in J_i}\bar r^j_iy_j\big) \Big) .
\label{p12}
\end{align}
Since~\eqref{p11} is a Presburger
formula of $(\Zpm,0,+,\leq)$,
so does~\eqref{p12}, so Claim C
is proved. 
Therefore, every rational
set of $((\Zpm)^n,+)$ is definable 
by a Presburger formula of $(\Zpm,0,+,\leq)$.

Let us now show that every subset of $(\Z\cup\{\pm\infty\})^n$
definable by a Presburger formula of $(\Zpm,0,+,\leq)$
is a rational set of the commutative monoid $((\Zpm)^n,+)$.
As the family of rational sets of
$((\Zpm)^n,+)$
is closed with respect to union, intersection and complementation
(see \cite{eilenberg69}) 
and as the projection of a rational set of $((\Zpm)^n,+)$
is a rational set, it is enough to show that 
for all nonnegative integers, $r_i$, $k_i$, $1\leq i \leq n$, 
the set $S$ of solutions of 
\begin{eqnarray}
\sum_{i=1}^n k_ix_i\leq \sum_{j=1}^n r_jx_j \label{fatomic}
\end{eqnarray}
is a rational set of $((\Zpm)^n,+)$.
To see this, consider the map
$\kappa:\Zpm\to\pattern$,
which fixes $-\infty$ and $+\infty$,
and sends any finite number to $0$.
We extend $\kappa$ to a map
$(\Zpm)^n\to\pattern^n$ by
making $\kappa$ act on each entry.
We shall call $\kappa(x)$ the 
{\em pattern} of $x\in (\Zpm)^n$.
To show that $S$ is rational,
it is enough to prove that for every
$p\in \pattern^n$,
the set of solutions with pattern $p$,
$S_p=S\cap \kappa^{-1}(p)$,
is rational. Let $I(p)=\set{1\leq i\leq n}{p_i=\pm\infty}$,
and let $J(p)$ denote the complement of $I(p)$
in $\{1,\ldots,n\}$. When $x$
has pattern $p$,~\eqref{fatomic} can be rewritten
as 
\begin{eqnarray}
a+ \sum_{i\in J(p)} k_ix_i\leq b+ \sum_{j\in J(p)} r_jx_j 
\enspace,
\label{fatomic-new}
\end{eqnarray}
where $a= \sum_{i\in I(p)} k_i p_i$
and $b= \sum_{j\in I(p)} r_j p_j$. 
Note that $a$ and $b$ can take values only in $\pattern$
($a=0$ if $k_i=0$ for all $i\in I(p)$,
and $b=0$ if $r_j=0$ for all $j\in I(p)$,
due to the convention $0\times x_i=0$
which is implied by the convention $x_i^0=e$
that we made when writing~\eqref{p1}).
Note also that an empty sum is equal to $0$,
by convention. Only the following cases can occur.
\newcommand{\mycase}[1]{\medskip\noindent{\em Case #1.}\/}

\mycase 1
$a=-\infty$. Then, ~\eqref{fatomic-new} is
automatically satisfied, and we conclude that  
$S_p=\kappa^{-1}(p)=\set{x\in (\Zpm)^n}{x_i=p_i\mrm{ for }i\in I(p)\mrm{ and }x_i\in \Z\mrm{ for }i\in J(p)}$.
The set $S_p$, which is a Cartesian product
of rational sets of commutative monoids,
is rational. (Indeed, $S_p$ is the Cartesian product
of one element subsets of $(\Zpm,+)$,
which are obviously rational,
and of copies of $\Z=\{-1,1\}^*$
which is a rational subset of the monoid $(\Zpm,+)$.)

\mycase 2  $a\neq-\infty$. We split this case
into subcases:

\mycase{2.1} $a\neq-\infty$ and $b=-\infty$.
Then, $S_p=\emptyset$.

\mycase{2.2}
$a\neq-\infty$ and $b=+\infty$. Then,
$S_p=\kappa^{-1}(p)$, and we proved in Case~1 that $S_p$ is rational.

\mycase{2.3}
$a\neq-\infty$ and $b=0$. We again split
this subcase.

\mycase{2.3.1}
$a=+\infty$ and $b=0$. Then, $S_p=\emptyset$.

\mycase{2.3.2}
$a=0$ and $b=0$.
Then, $S_p=\set{x\in (\Zpm)^n}{x_i=p_i\mrm{ for }i\in I(p),
\; x_i\in \Z\mrm{ for }i\in J(p), \mrm{ and } \sum_{i\in J(p)}k_ix_i\leq \sum_{j\in J(p)} r_jx_j }$. 
By the classical result of~\cite{GinSpa}, the set
of {\em finite} integer solutions of an equation
of the form~\eqref{p11} is rational,
therefore,
$T=\set{x\in \Z^{J(p)}}{\sum_{i\in J(p)}k_ix_i\leq \sum_{j\in J(p)} r_jx_j}$ is a rational subset of
 $\Z^{J(p)}$.
Since $S_p$ is the Cartesian product of $T$ by
one element sets, $S_p$ is rational.

Thus, the set $S=\bigcup_{p\in \pattern^n}S_p$
of solutions of (\ref{fatomic})
is a rational set of  $((\Zpm)^n,+)$.
\end{proof}

\begin{example}
The idempotent semirings $\rmax$ and $\qmax=(\Q\cup\{-\infty\},\max,+)$
do not have the Presburger property. 
As a first counter example, 
consider the set $D_1=\set{x\in \cS}{x\geq 0}$,
where $\cS=\rmax$ or $\qmax$.
This set is defined by a Presburger formula,
but is not rational. 
Indeed, it follows from the representation
of rational sets as semilinear
sets that the set of non-zero elements
of any rational set contained
in $D_1$ has a minimal element, whereas
$D_1$ does not have this property.
To give a second counter example, note that any set definable by 
a Presburger formula of $(\Rm,0,+,\leq)$ or
$(\Q\cup\{-\infty\},0,+,\leq)$ must be stable by multiplication by a 
positive constant. Therefore, the rational set $D_2=\{1\}$
cannot be defined by a Presburger formula.
Another example of idempotent semiring which
does not have the Presburger
property is $\cS=((\Z\cup\{-\infty\})^2,\max,+)$,
where $\max$ denotes the law $(\Z\cup\{-\infty\})^2\times
(\Z\cup\{-\infty\})^2\to (\Z\cup\{-\infty\})^2$
which does entrywise max. In this semiring
the set $\{(1,0)\}$
cannot be defined
by a Presburger formula (for symmetry reasons, there is no way
to distinguish $(1,0)$ from $(0,1)$ using Presburger
formulas). 
\end{example}

\section{Closure Properties of Rational Semimodules}
\label{sec-closure}
Let us recall some definitions.
A (left) \NEW{semimodule} over 
a semiring $(\cS,\oplus,\otimes,\zero_{\cS},\unit)$
is a commutative monoid $(\cX,\oplus,\zero_{\cX})$,
equipped with a map $\cS\times \cX\to \cX$, $(\lambda,x)
\to \lambda x$ (left action), that satisfies
\begin{subequations}
\begin{gather}
(\lambda \otimes \mu) x= \lambda (\mu x)\enspace,\label{e-lefta}\\
\lambda(x\oplus y) = \lambda x\oplus \lambda y,\;
(\lambda \oplus \mu)x = \lambda x\oplus
\mu x \enspace,\label{e-dist0}\\
\zero_{\cS} x = \zero_{\cX},\;
\lambda \zero_{\cX} = \zero_{\cX} ,\;
\unit x =x \enspace,\label{e-unit}
\end{gather}
\end{subequations}
for all $x,y\in \cX$, $\lambda,\mu\in \cS$. 
In the sequel, we will denote 
by $\zero$ both the zero element $\zero_{\cS}$ of $\cS$
and the zero element $\zero_{\cX}$ of $\cX$, when there will be no risk
of confusion. We will also use concatenation
to denote the product of $\cS$, so that~\eqref{e-lefta}
will be rewritten as $(\lambda \mu) x= \lambda (\mu x)$.
When $(\cS,\oplus)$ is idempotent, 
$(\cX,\oplus)$ is idempotent
(indeed, it follows from~\eqref{e-dist0}
and~\eqref{e-unit} that $x=\unit x=
(\unit\oplus\unit)x=\unit x\oplus \unit x=x\oplus x$).
A \NEW{subsemimodule} of $\cX$ is 
a subset $\cZ\subset \cX$ such that
$\lambda x \oplus \mu y \in \cZ$,
for all $x,y\in \cZ$ and $\lambda,\mu\in \cS$.
We will consider subsemimodules
of the \NEW{free semimodule} $\cS^n$, which
is the set of $n$-dimensional vectors over $\cS$,
equipped with the internal law $(x\oplus y)_i=x_i\oplus y_i$
and the left action $(\lambda x)_i=\lambda \otimes x_i$.
If $G\subset {\cS}^n$, 
we will denote by 
$\vect G $ the subsemimodule of ${\cS}^{n}$ generated by $G$, i.e. the set of all $x\in {\cS}^{n}$ for which there exists a finite number of elements $u^1,\ldots ,u^k$ of $G$ and a finite number of scalars $\lambda_1,\ldots ,\lambda_k\in {\cS}$, such that $x=\bigoplus_{i=1}^k\lambda_i u^i$.
\begin{definition}[Rational semimodules]\label{def-rat}
A subsemimodule ${\cX}\subset {\cS}^{n}$ is \NEW{rational}
if it has a generating family which is a rational subset
of the monoid $({\cS}^{n}, \otimes)$. 
\end{definition}
We now show that rational semimodules
are closed under natural algebraic
operations.
We begin by a simple general
property.
\begin{theorem}\label{theosimple}
Let ${\cS}$ be an arbitrary semiring. Let ${\cX}$, ${\cY}\subset {\cS}^{n}$ and ${\cZ}\subset {\cS}^{p}$ be rational semimodules. Then ${\cX}\oplus {\cY}$ and ${\cX}\times {\cZ}$ are rational semimodules.
\end{theorem}
\begin{proof}
Let $\cX=\vect B$, $\cY=\vect C$ and $\cZ=\vect D$,
where $B,C\subset\cS^n$ and $D\subset\cS^p$ are rational sets.
As $\cX\oplus\cY=\vect( B\cup C)$, it follows that
$\cX\oplus \cY$ is a rational semimodule because 
$B\cup C$ is a rational set of the monoid
$(\cS^n,\otimes )$.
Let us denote by $\zero_{k}$, for $k\in \N$,
the neutral element for addition in the commutative monoid
$(\cS^k,\oplus)$ and let us define the sets
$$\tilde{B} =\left\{ \begin{pmatrix}x \\ \zero _{p}\end{pmatrix}
\mid x\in B\right\} \mrm{ and } \hat{D} =\left\{ \begin{pmatrix}
\zero _{n} \\ z \end{pmatrix}\mid z\in D\right\} .$$
Since $B$ and $D$ are rational sets, $\tilde{B} $ and $\hat{D} $ are rational sets of the monoid $(\cS^{n+p},\otimes)$.
Then, as $\cX\times\cY=\vect(\tilde{B} \cup \hat{D})$,
it follows that ${\cX} \times {\cY}$ is a rational semimodule.
\end{proof}

We will need the following analogue of Caratheodory theorem, which
was already stated in~\cite{gaubert98n}.
(The classical Caratheodory theorem,
see e.g.~\cite[Cor.~7.1i]{schrijver},
states that if a vector $x$ in $\R^n$ is
a positive linear combination of vectors
of a finite subset $G\subset \R^n$, $x$ can
be written as a positive linear combinations of at most
$n$ vectors of $G$). 
\begin{proposition}[Max-Plus Caratheodory Theorem]\label{pc} 
\sloppy Let $\cS$ be an idempotent semiring whose natural order is a total
order.
If $G\subset\cS^{n}$, and if $x\in ${\rm span}$\left( G\right)$, then there is a subset $B$ of $G$, of cardinality at most $n$, such that $x\in ${\rm span}$\left( B\right)$.
\end{proposition}
\begin{proof}
If $x\in \vect G$, 
we can find  $u^{1},\ldots , u^{k} \in G, \lambda_1,\ldots,\lambda_k
\in \cS$ such that $x=\oplus _{1\leq i\leq k}\lambda_i u^i$. 
For all $1\leq j\leq n$,
we define $I(j) = \set i{1\leq i\leq k, \;x_j=\lambda_i u_j^i}$.
Since the natural order of $\cS$ is a total order,
all the $I(j)$ are non-empty.
Choosing exactly one element $i_j$ in $I(j)$,
we obtain a family  $u^{i_1},\ldots ,u^{i_n}$ such that
$x=\lambda _{i_1}u^{i_1}\oplus\cdots\oplus\lambda_{i_n}u^{i_n}$
(It may be the case that $i_j=i_k$ for some $j\neq  k$. 
In such cases the cardinality of $B$ is less than $n$.)
\end{proof}
   
\begin{theorem} \label{rateqratset}
Let ${\cS}$ be an
idempotent commutative semiring which
satisfies the Presburger property
and whose natural order is a total order.
Then, for a subset $\cX$ of $\cS^n$
the following statements are equivalent.
\begin{enumerate}
\item\label{it1} $\cX$ is a rational semimodule.
\item\label{it2}
 $\cX$ is a semimodule and a rational subset of the monoid
$(\cS^n,\otimes)$.
\end{enumerate}
\end{theorem} 
\begin{proof}
Trivially, (\ref{it2}) implies (\ref{it1})
because $\cX$ is generated by $\cX$.
Suppose now that $\cX$ is a rational semimodule and let
$G\subset\cS^n$ be a rational set such that
$\cX=\vect G$. Let $P$ be a formula of the first-order logic of
$(\cS,e,\otimes,\preceq)$, that defines $G$.
The Max-Plus Caratheodory Theorem implies that:
$x\in \cX$ if and only if 
\begin{align*}
&(\exists u^1\in \cS^n) ,\ldots,(\exists u^n\in\cS^n),
(\exists\lambda_1\in\cS) ,\ldots,(\exists\lambda_n\in\cS),
\\
&\qquad \qquad\qquad \Big( P(u^1)\wedge \cdots \wedge P(u^n)\wedge x=\bigoplus_{i=1}^{n} {\lambda }_i u^i\Big) \enspace.
\end{align*}
Since the last formula belongs to the first-order logic of 
$(\cS,e,\otimes,\preceq)$, 
we obtain that $\cX$ is a rational set of the commutative monoid
$(\cS^n,\otimes)$
\end{proof}
If ${\cX}$ and ${\cY}$ are two semimodules over $\cS$, we denote by
 $\Hom(\cX,\cY) $ the set of \NEW{linear maps}, i.e.,
of semimodule morphisms, from $\cX$ to $\cY$.
A linear map $\cS^n\to\cS^p$
can be represented uniquely in matrix form,
$x\mapsto Ax, \; (Ax)_i=\bigoplus _{1\leq j\leq n} A_{ij}x_j$,
where $A=(A_{ij}) \in\cS^{p\times n}$.

\begin{theorem}[Closure theorem]\label{closureth}
Let $\cS$ be an idempotent commutative semiring which satisfies
the Presburger property and  whose natural order is a total order.
Let $\cX$, $\cY\subset\cS^n$,
$\cZ\subset\cS^p$, 
$\cG\subset\cS^{n+p}$ and
$\cW\subset(\cS^n)^2$ be rational semimodules,
and let $A\in \Hom(\cS^n,\cS^p)$. 
Then the following sets all are rational semimodules.
\begin{enumerate}
\item\label{ip1} $\cX\cap\cY$,
\item\label{ip2}
$\cX\cG=\set{v\in\cS^p}{\exists x\in\cX,(x,v)\in\cG}$
and $\cG\cZ=\set{u\in\cS^n}{\exists z\in\cZ,(u,z)\in\cG}$,
\item\label{ip3}
$A\cX=\set{Ax}{x\in \cX}$,
\item\label{ip4}
$A^{-1}\cZ=\set{x\in \cS^n}{Ax\in \cZ}$,
\item\label{ip5}
$\cX\ominus \cY=
\set{u\in\cS^n}{\exists y\in\cY, u\oplus y\in \cX}$,
\item\label{ip6}
$\cW^\perp=
\set{x\in \cS^n}{a\cdot x=b\cdot x,\forall (a,b)\in \cW} $,
where $a\cdot x=\bigoplus_{1\leq i\leq n}a_i x_i$,
\item\label{ip7}
$\cX^\top=
\set{(a,b)\in (\cS^n)^2}{a\cdot x=b\cdot x,\forall x\in \cX}$.
\end{enumerate}
\end{theorem}
\begin{proof}
\ref{ip1}. Since $\cX$ and $\cY$ are rational semimodules, 
we know that they are rational sets 
(by Theorem~\ref{rateqratset}). 
As the intersection of rational sets of a commutative monoid
is a rational set (see \cite{eilenberg69}), 
we have that $\cX\cap\cY$ is a rational
set and therefore a rational semimodule.

\ref{ip2}. By symmetry, we only consider the case of $\cX\cG$.
Since $\cG\subset\cS^{n+p}$ and  $\cX\subset \cS^n$
are rational semimodules, we know that they are rational sets
of $(\cS^{n+p},\otimes)$ and $(\cS^n,\otimes)$
respectively (by Theorem~\ref{rateqratset}). 
Let $P$ and $Q$ be Presburger formulas of
$(\cS,e,\otimes,\preceq)$,
defining $\cG$ and $\cX$ respectively. 
Then as $$v\in\cX\cG \mrm{ if and only if }
(\exists x\in\cS^n) ( Q(x)\wedge P(x,v)),$$
it follows that $\cX\cG$,
which is defined by a Presburger formula of
$(\cS,e,\otimes ,\preceq)$,
is a rational set. By Theorem~\ref{rateqratset},
it is a rational semimodule.

\ref{ip3}. Let us define $\cG=\set{(u,Au)}{u\in \cS^n}$.
Since $\cG$ is finitely generated,
$\cG$ is a rational semimodule,
and $A\cX=\cX\cG$ is a rational semimodule.

\ref{ip4}. Taking again $\cG=\set{(u,Au)}{u\in \cS^n}$,
we have $A^{-1}\cZ=\cG\cZ$. 
Hence, $A^{-1}\cZ$ is a rational semimodule.

\ref{ip5}.
Let us define ${\cG}=\set{(u,y,x) \in  (\cS^n)^3}{x=u\oplus y}$.
Then $\cG$, which is the direct image of $(\cS^n)^2$
by a linear map,
is finitely generated, and a fortiori,
rational. 
Thus, $\cX\ominus\cY = 
\set{u\in\cS^n}{\exists x \in \cX,\exists y\in \cY, x=u\oplus y }
=
\set{u\in\cS^n}{\exists x\in\cX,\exists y\in\cY,(u,y,x)\in \cG}
= \cG(\cY\times \cX)$ is a rational
semimodule.

\ref{ip6}. As $\cW\subset(\cS^n)^2$ is a rational semimodule, 
we know (by Theorem~\ref{rateqratset}) that it is a rational set.
Let $P(u_1,\ldots,u_n,v_1,\ldots,v_n)$ be a Presburger
formula of $(\cS,e,\otimes,\preceq)$ defining $\cW$.
Then $(x_1,\ldots,x_n)\in \cW^\perp$ if and only if
\begin{align*}
&\forall u_1\in\cS,\ldots ,\forall u_n\in\cS,\forall v_1\in \cS,\ldots ,\forall v_n\in\cS \\
&\qquad 
\big( P( u_1,\ldots ,u_n,v_1,\ldots ,v_n) \Rightarrow
 \bigoplus_{i=1}^n x_iu_i=\bigoplus_{j=1}^n x_jv_j \big) .
\end{align*}
Since this is a Presburger formula of $(\cS,e,\otimes,\preceq)$,
it follows that $\cW^{\perp}$ is a rational set of
$(\cS^n,\otimes)$, and also, by Theorem~\ref{rateqratset},
a rational semimodule.

\ref{ip7}. 
Let $P(x_1,\ldots ,x_n)$ be a Presburger formula of 
$(\cS,e,\otimes,\preceq) $ defining $\cX$. Then we have that 
$(u,v) \in \cX^{\top}$ if and only if
\[
(\forall x_1\in\cS) ,\ldots ,(\forall x_n\in \cS),
(P(x_1,\ldots ,x_n) \Rightarrow
\bigoplus_{i=1}^n u_i x_i =\bigoplus_{j=1}^n v_i x_i) \enspace. 
\]
Arguing as in Statement~\ref{ip6},
we conclude that $\cX^{\top}$
is a rational semimodule.
\end{proof}
\begin{remark}
A motivation for considering the operations
$\ominus$ and $\cZ\to A^{-1}\cZ$ comes
from $(A,B)$ invariant spaces (see~\cite{wonham}).
If one consider the dynamical system
\[
x(k)=Ax(k-1)\oplus Bu(k),
\]
where $A\in(\zmax)^{n\times n}$, 
$B\in(\zmax)^{n\times p}$,
$x(k)\in (\zmax) ^n$,
and $u(k)\in(\zmax)^p$,
the set of $x(0)$ for which
there exists a control $u(1)$
such that $x(1)$ belongs to 
a prescribed semimodule $\cX$
is $A^{-1}(\cX\ominus \cB)$,
where $\cB$ denotes the semimodule
generated by the columns of $B$.
Max-plus $(A,B)$-invariant spaces are further studied
in~\cite{katz}.
\end{remark}
We shall say that a vector $v$ of a semimodule $\cX\subset \cS^n$
is \NEW{extremal} if $v\not\in \vect(\cX \setminus\vect \{v\})$.
We denote by $\mrm{Ext}(\cX)$ the set
of extremal points of $\cX$.
The interest in extremal points stems from
a theorem due to Moller~\cite{moller88} and Wagneur~\cite{wagneur91},
which states that a finitely generated subsemimodule
of $(\rmax)^n$ is generated by its extremal vectors. 
\begin{theorem}\label{extrat}
Let $\cS$ be an idempotent commutative semiring which
satisfies the Presburger property and whose natural order is
a total order. 
If $\cX\subset\cS^n$ is a rational semimodule, then 
$\Ext(\cX)$
 is a rational set of the monoid $(\cS^n,\otimes)$.
\end{theorem}
\begin{proof}
Let $P$ be a Presburger formula of
 $(\cS,e,\otimes ,\preceq)$
defining $\cX$.
The max-plus Caratheodory theorem shows
that $v\in \vect(\cX\setminus\vect\{v\})$
is equivalent to
\begin{align*}
&(\exists u^1\in \cS^n), \ldots ,(\exists u^n\in \cS^n),
( \exists \lambda_1\in \cS) , \ldots , 
(\exists \lambda_n\in \cS) 
 P(u^1)\wedge\cdots\wedge P(u^n)\\
&\quad \wedge (v=\bigoplus_{i=1}^{n} \lambda _i u^i )
\wedge
\neg \left( (\exists\alpha_1\in\cS)(u^1=\alpha_1v)\vee \cdots 
\vee(\exists \alpha_n\in\cS)(u^n=\alpha_nv)\right)
\enspace .
\end{align*}
Since this is a Presburger formula of $\left( {\cS}, e,\otimes ,\preceq \right) $, it follows that $\cX\setminus\Ext(\cX)$ is a rational set, 
and therefore $\Ext(\cX)$ is a rational set.
\end{proof}
\begin{remark}
We could prove Statement~\ref{ip3}
of Theorem~\ref{closureth} without using
Presburger's arithmetics, as follows. If $R$ is a rational
set that generates the semimodule $\cX$,
$A\cX$ is generated by the set $A(R)=\set{Ar}{r\in R}$. One can show
directly, using the fact that a max-plus linear map
is piecewise affine with integer slopes, that $A(R)$ is rational.
\end{remark}
\begin{remark}\label{rk-algo}
A difficulty, in looking for more
direct proofs of 
Statements~\ref{ip1},\ref{ip2},\ref{ip4}--\ref{ip7}
of Theorem~\ref{closureth},
is the relative absence of knowledge
of the minimal set of generators of
a semimodule defined by natural algebraic
operations. 
This difficulty persists even in the case of
finitely generated semimodules.
For instance, the only known algorithm 
(see~\cite{butkovicH}, ~\cite[III,1.1.4]{gaubert92a} 
or~\cite[Th.~8]{maxplus97})
to compute a generating family of the set of solutions
of the max-plus linear system $Ax=Bx$,
where $A,B$ are $n\times p$ matrices,
has an a priori doubly exponential execution time,
and tells little about the geometry of extremal points.
(However, the doubly exponential bound
is pessimistic,
the average case is better in practice, and finding
only one solution can be done more efficiently
by computing sub-fixed point of min-max functions,
see~\cite{walkup,maxplus97,bcg99} and~\cite{cras,gg}
for fixed point algorithms for min-max functions.)
\end{remark}
\section{Examples and Counter Examples}
\label{secAyE}
\subsection{Reachable and Observable Spaces of Max-Plus Linear Discrete Event Systems}\label{subsec-obs}
Let us consider the max-plus linear system:
\begin{subequations}
\begin{align}
x(k)&=Ax(k-1)\oplus Bu(k),\label{mata}\\
y(k)&=Cx(k) ,\\ 
x(0)&=\xi \enspace .
\end{align}\label{mpsystem}\end{subequations}
where $A\in (\zmax)^{n\times n}$,
$B\in(\zmax)^{n\times p}$,
$C\in(\zmax)^{q\times n}$,
$\xi\in (\zmax) ^n$,
and $u(k)\in(\zmax)^p$, $k=1,2,\ldots$ is a sequence of control vectors.
We call {\em reachable space} in time $k$,
and denote by $\cR_k$, the set of states $x(k)$
reachable from the initial state $x(0)=\zero$.
We also define the \NEW{reachable space} in arbitrary time,
$\cR_{\omega}$, which is the union of the $\cR_k$.
(We shall sometimes write $\cR_k(A,B)$
or $\cR_{\omega}(A,B)$ to emphasize the
dependence in $A,B$.)
Introducing the \NEW{reachability matrices}
\[\Rmat_k=(B,AB,\ldots, A^{k-1}B),
\qquad \Rmat_{\omega}=(B,AB,A^2B\ldots)
\enspace,
\]
we characterize $\cR_{k}$ (resp. $\cR_{\omega}$)
as the semimodule generated
by the columns of the matrix
$\Rmat_k$ (resp. $\Rmat_{\omega}$).
Identifying matrices with operators,
we will write $\cR_k=\im\Rmat_k$,
and $\cR_{\omega}=\im\Rmat_{\omega}$.

The definition of rational semimodules
is motivated by the following result:
\begin{theorem}\label{reachableth}
Reachable spaces are rational semimodules, i.e.
if $A\in (\zmax)^{n\times n}$ and $B\in(\zmax)^{n\times p}$, 
then $\cR_{\omega}=\im \Rmat_{\omega}$ is a rational semimodule.
\end{theorem}
The proof of Theorem~\ref{reachableth}
relies on the following cyclicity theorem
for reducible max-plus matrices,
which is taken from~\cite[VI,1.1.10]{gaubert92a}.
\begin{theorem}\label{cycli}
Let $A\in (\rmax)^{n\times n}$. 
There are positive integers $c,N$, such that for all
$1\leq i,j\leq n$, 
there are scalars $\lambda_0,\ldots, \lambda_{c-1}$
(depending on $i,j$) such that
for all $0\leq l\leq c-1$,
\begin{align}
\forall n\geq N,\qquad 
(A^{nc+l+c})_{ij}=\lambda_l (A^{nc+l})_{ij} \enspace .
\label{e-cycli}
\end{align}
\end{theorem}
This cyclicity theorem follows readily from
the characterization of max-plus rational series
in one variable as merge of ultimately
rational series,
see~\cite{moller88}, \cite[VI,1.1.8]{gaubert92a}, \cite{krob93a},
and the discussions in~\cite{gaubert94d,maxplus97}.

\begin{proof}[Proof of Theorem~\ref{reachableth}]
Theorem~\ref{cycli} implies
that $\{A^0,A^1,A^2,\ldots\}$ is a rational
subset of $((\Zm)^{n\times n},+)$, and
therefore, $\cX=\vect\{A^0,A^1,A^2,\ldots\}$
is a rational subsemimodule of $(\zmax)^{n\times n}$.
Since $\cR_{\omega}(A,B)$
is the sum of the reachability 
spaces $\cR_{\omega}(A,B_{\cdot,i})$
associated to the differents columns
$B_{\cdot,i}$ of $B$, for $1\leq i\leq p$, and
since the sum of rational semimodules
is rational (cf. Theorem~\ref{theosimple}),
it is enough to consider
the case when $B$ has only one column.
Then, $\cR_{\omega}$ is the direct image of $\cX $ by
the linear map $(\zmax)^{n\times n}\to (\zmax)^{n\times 1},
X\mapsto XB$,  
and it follows from Statement~\ref{ip3} of Theorem~\ref{closureth}
that $\cR_{\omega}$ is rational.\end{proof}
Let $\xi,\xi'\in (\zmax)^n$,
and consider two trajectories 
of the dynamical system~\eqref{mpsystem}, 
\[ \{(x(k),y(k))\}_{k\geq 0},
\qquad\mrm{and}\qquad\{(x'(k),y'(k))\}_{k\geq 0}
\enspace,
\]
corresponding
to the initial conditions $x(0)=\xi$, $x'(0)=\xi'$,
the zero control $u(k)\equiv \zero$
being applied in both cases. 
We call \NEW{observable congruence} in time $k\geq 1$,
and denote by $\cO_k$,
the congruence over $(\zmax)^n$ defined by
\[
(\xi, \xi')\in \cO_k 
\iff y(l)=y'(l),\;\forall 0\leq l\leq k-1 \enspace ,
\]
and the \NEW{observable congruence} (in arbitrary time)
$\cO_{\omega}$ is defined
as the intersection of the congruences $\cO_k$,
$k\geq 1$. By congruence, we mean
an equivalence relation on $(\zmax)^n$ compatible
with the semimodule structure of $(\zmax)^n$.
In particular, $\cO_k$ and $\cO_{\omega}$
are subsemimodules of $((\zmax)^n)^2$.
Introducing the \NEW{observability matrices}
\[
\Omat_k= 
\begin{pmatrix}
C \\
CA \\ 
\vdots \\
CA^{k-1}
\end{pmatrix}
,\qquad
\Omat_{\omega}= 
\begin{pmatrix}
C \\
CA \\ 
CA^2\\
\vdots
\end{pmatrix}\enspace,
\]
we characterize $\cO_k$ (resp. $\cO_{\omega}$) as the right kernel
$\ker \Omat_k$ (resp. $\ker \Omat_{\omega}$)
of $\Omat_k$ (resp. $\Omat_{\omega}$), that is:
\[
(\xi,\xi')\in \cO_k \iff \Omat_k\xi=
\Omat_k\xi',\qquad
(\xi,\xi')\in \cO_{\omega} \iff \Omat_{\omega}\xi=
\Omat_{\omega}\xi'\enspace .
\]
See~\cite{ccggq99} for more background on max-plus
reachability spaces and observable congruences.
We have the following dual version of Theorem~\ref{reachableth}:
\begin{theorem}\label{observath}
Observable congruences are rational,
i.e. if $A\in (\zmax)^{n\times n}$, 
$C\in (\zmax)^{q\times n}$,
then $\cO_{\omega}=\ker \Omat_{\omega}$ is a rational subsemimodule
of $((\zmax)^n)^2$. 
\end{theorem}
\begin{proof}
By Theorem~\ref{reachableth},
the semimodule $\cZ$ 
generated by the rows of the observability
matrix $\Omat_{\omega}$,
which can be identified to the reachable space
$\cR_{\omega}(A^T,C^T)$, is rational.
Since $\cO_{\omega}=\cZ^{\top}$,
Statement~\ref{ip7} of Theorem~\ref{closureth}
shows that $\cO_{\omega}$ is rational.
\end{proof}
\subsection{Example of reachable space and observable congruence}\label{subsec-display}
\begin{figure}
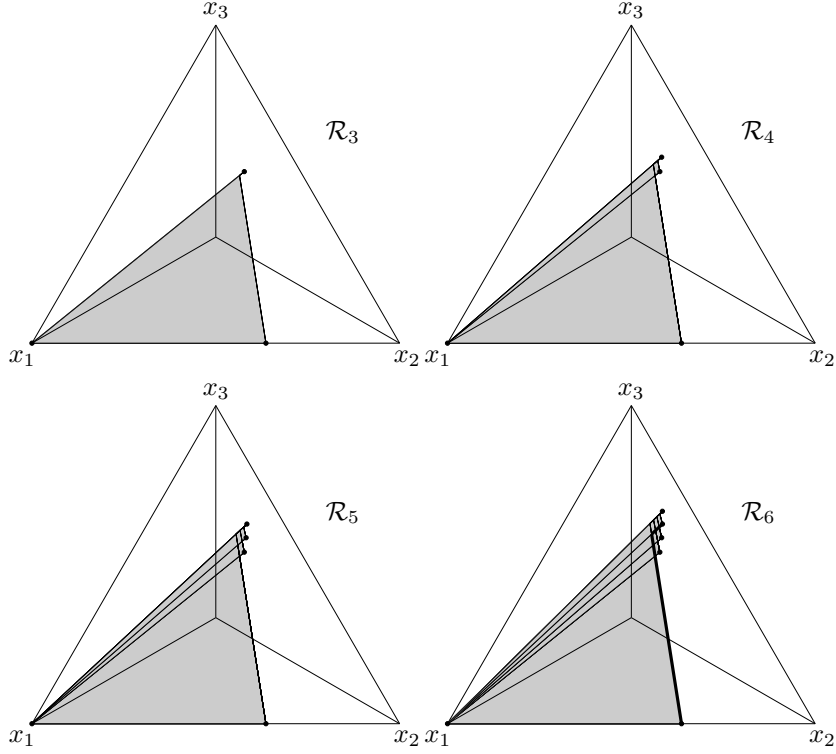

\begin{center}
\begin{tabular}[t]{cc}
\input fig1 &
\input fig2 \\
\input fig3 &
\input fig4
\end{tabular}
\end{center}
\caption{Exponential representation of the reachable spaces $\cR_3,\cR_4,\cR_5,\cR_6$ for the $(A,B)$ pair given by~\eqref{exab}}
\label{fig-semi}
\end{figure}

Consider
\begin{align}
A=\begin{pmatrix}1 & -\infty & -\infty \cr 
5& 2 & -\infty \cr
-\infty & 6 & 3 \end{pmatrix},
\qquad 
B=\begin{pmatrix}0\cr -\infty\cr -\infty\end{pmatrix}
\enspace .
\label{exab}
\end{align}
Then $\cR_{\omega}=\im\Rmat_{\omega}$ where
\begin{align}
\Rmat_{\omega}=\begin{pmatrix}
0       & 1       & 2  & 3  & 4  & 5  & 6  & \cdots \cr
-\infty & 5       & 7  & 9  & 11 & 13 & 15 & \cdots \cr
-\infty & -\infty & 11 & 14 & 17 & 20 & 23 & \cdots \end{pmatrix}\enspace .
\label{exr}
\end{align}
Obviously $\cR_{\omega}$ is a rational semimodule because
the set of columns of $\Rmat_{\omega}$ can be written
as $U\cup (\{v\}+\{w\}^*)$, with
\begin{align}\label{gen}
U= \left\{
\begin{pmatrix}
0\cr
-\infty\cr
-\infty \cr
\end{pmatrix},
\begin{pmatrix}
1\cr      
5\cr
-\infty 
\end{pmatrix}
\right\}
,
\qquad
v=\begin{pmatrix}
2\cr 7  \cr 11
\end{pmatrix}
,\qquad
w=
\begin{pmatrix}
1\cr 
2\cr
3\end{pmatrix}
\end{align}
The semimodules $\cR_3,\cR_4,\cR_5,\cR_6$ are shown
on Figure~\ref{fig-semi}.

To represent semimodules
which contain vectors
with $-\infty$ coordinates,
we use the following projection.
Let $\beta$ denote a positive parameter,
and let us chose a triangle in the plane.
We represent a point $x\in (\R\cup\{-\infty\} )^3$
by the point $\pi(x)$ of the triangle, which is
the barycenter of the vertices of the triangle with respective
weights $\exp(\beta x_1)$, $\exp(\beta x_2)$,
$\exp(\beta x_3)$. We shall refer to this
projection as the \NEW{exponential projection}
in the sequel. The exponential projection has
the property that if two points $x$ and $y$
are proportional in the max-plus sense,
that is, if $x_i=\lambda+y_i$ for some $\lambda\in \R$,
then, $\pi(x)=\pi(y)$. Therefore, representing the image
of a semimodule $\cX\subset(\rmax)^3$
(or $\cX\subset(\zmax)^3$) by $\pi$ gives
a faithful image of $\cX$. Such drawings
represent in fact the max-plus two-dimensional
\NEW{projective space}, which is the quotient
of $(\rmax)^3$ by the parallelism relation $x\sim y
\iff x=\lambda+y$ for some $\lambda\in \R$.
The max-plus projective space appeared in
the work of several authors, 
see~\cite{kolokoltsov,gaubert95c,mairesse95a,gaubert98n}.

In Figure~\ref{fig-semi}, the generators of the semimodules
$\cR_3,\cR_4,\cR_5,\cR_6$,
that is, the columns of the matrices 
$\Rmat_3,\Rmat_4,\Rmat_5,\Rmat_6$,
are represented by bold points. For any two generators,
we have represented the max-plus plane
generated by these
two generators
(we call \NEW{plane} a semimodule
generated by two nonproportional vectors).
The projection $\pi$ sends in general
a plane to a broken segment. For instance,
the bold broken segment
on the fourth picture of Figure~\ref{fig-semi}
represents the max-plus plane generated by the second and
fifth columns of $\Rmat_6$.

It should be graphically clear from
Figure~\ref{fig-semi} that the 
generators are extremal, that the
semimodules $\cR_k, k=0,1,2,3,\ldots$ form
an infinite ascending sequence (this illustrates
the fact that the semimodule $(\zmax)^3$
is not Noetherian),
and that ${\cR_{\omega}}$ is not finitely generated. 
One can check mechanically all these facts by appealing
to residuation theory, which allows
us to compute the extremal vectors
of semimodules, see~\cite{bcoq},\cite{butkovip94} and~\cite{maxplus97}
for more details. Let us also mention that the computations
of this example have been checked
using the max-plus toolbox of scilab,
see~\cite{toolbox}.

We can visualize, on the drawings of Figure~\ref{fig-semi},
both the $\rmax$ semimodule and the $\zmax$ semimodule
generated by the columns of the matrices $\Rmat_k$.
The gray zone represent an $\rmax$ semimodule.
The corresponding $\zmax$ semimodule is
an ``integer lattice''
inscribed in the real semimodule,
that for readability of the figure, we
do not have represented. 

\begin{figure}
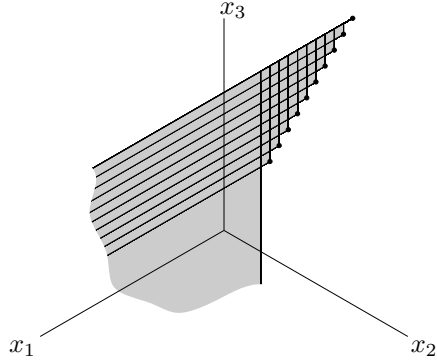

\begin{center}
\input fig5
\end{center}
\caption{Orthogonal projection of the reachable space $\cR_{12}$ for the $(A,B)$ pair of ~\eqref{exab}}
\label{figorth}
\end{figure}

To see graphically that the semimodule
$\cR_\omega$ is rational, it is convenient
to use another representation, in
which every {\em finite} point of $\cR_{\omega}$ is projected
orthogonally to the main diagonal of $\R^3$:
again, two vectors $x,y\in \R^3$ which
are proportional in the max-plus sense,
are sent to the same point.
Using this projection,
the semimodule $\cR_{12}$ is represented
on Figure~\ref{figorth}. The rationality
of $\cR_{\omega}$ can be visualized on this figure:
the set of finite generators of $\cR_{\omega}$,
which is given by $\{v\}+\{w\}^*$,
where $v,w$ are as in~\eqref{gen},
is precisely the discrete half line of bold points.

Let us now represent an observability congruence. We consider the 
transposed dynamical system with new observation matrix $C=B^T$ and new
dynamics $A^T$. Then, the observability matrix is 
$\Omat_\omega(A^T,B^T)=(\Rmat_\omega(A,B))^T$, that is, the transpose of
the matrix computed in~\eqref{exr}. The corresponding observable
congruence $\mathcal{O}_{\omega }$ is depicted in Figure~\ref{fig-dual}, 
using the technique of~\cite[\S~4.3]{ccggq99}:
\begin{figure}
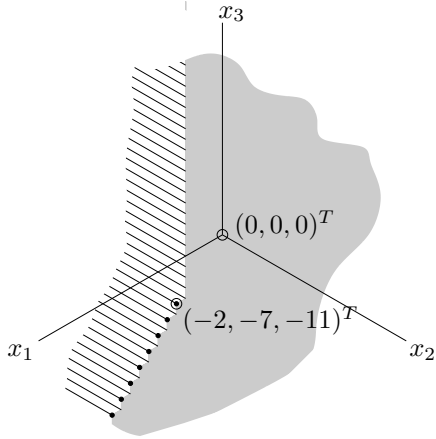

\begin{center}
\input fig6
\end{center}
\caption{Orthogonal projection of the observable congruence associated to $(A^T,B^T)$}
\label{fig-dual}
\end{figure}
(We only give the orthogonal projection here,
an exponential projection of another congruence
will be given later on, in Figure~\ref{fig-obs}.)
We know by  
Theorem~\ref{observath} that this is a rational congruence. The gray
region of Figure~\ref{fig-dual} represents the semimodule over the
min-plus semiring
$\zmin=(\Zp,\min,+)$ 
generated by the opposite of the rows of $\Omat_\omega$ (the
min-plus generators are represented by bold points): we can
derive from~\cite{ccggq99} that equivalence classes for
$\cO_\omega$ of points of the interior of this semimodule are
singletons. Let us check this elementarily 
for the point $\xi=(0,0,0)^T$ (indicated
by one of the two circles on the figure). We have:
$\Omat_\omega\xi=(0,5,11,14,17,\ldots )^T$. If 
$\Omat_\omega\xi=\Omat_\omega\xi'$ , from 
$(\Omat_\omega\xi')_1=(\Omat_\omega\xi)_1=0$  
it follows that $\xi_1'=0$. Then 
$(\Omat_\omega\xi')_2=\max (\xi_1'+1,
\xi_2'+5)=(\Omat_\omega\xi)_2=5$
implies that $\xi_2'=0$. Finally, from $(\Omat_{\omega
}\xi')_3=\max (\xi_1'+2,\xi_2'+7,\xi
_3'+11)=(\Omat_\omega\xi)_3=11$ it follows that $\xi
_3'=0$. Therefore $\Omat_\omega\xi=\Omat_{\omega
}\xi'\Rightarrow \xi'=\xi=(0,0,0)^T$. Other
equivalence classes are half-lines, as shown on the figure.
As an example 
let us compute the equivalence class of the point $\xi=(-2,-7,-11)^T$
(also indicated by a circle on the figure). 
We have that $\Omat_\omega\xi=(-2,-1,0,3,6,\ldots )^T$. If
$\Omat_\omega\xi=\Omat_\omega\xi'$, from 
$(\Omat_\omega\xi')_1=(\Omat_\omega\xi)_1=-2$ it
follows that $\xi_1'=-2$. Then $(\Omat_\omega\xi
')_2=\max (\xi_1'+1,\xi_2'+5)=(\Omat_\omega\xi)_2=-1$ implies that $\xi_2'\leq
-6$ and $(\Omat_\omega\xi')_3=\max (\xi_1'+2,\xi_2'+7,\xi_3'+11)=
(\Omat_\omega\xi
)_3=0$ implies that
$\xi_2'\leq -7$ and $\xi_3'\leq -11$. Finally, from
$(\Omat_\omega\xi')_4=\max (\xi_1'+3,\xi
_2'+9,\xi_3'+14)=(\Omat_\omega\xi)_4=3$
it follows that $\xi_3'=-11$. Now a straightforward
computation shows that any point $\xi'$ of the form 
$(-2,\alpha,-11)^T$, where $\alpha \leq -7$, satisfies
$\Omat_\omega\xi=\Omat_\omega\xi'$. Therefore the 
equivalence class of $\xi=(-2,-7,-11)^T$ is 
$\set{(-2,\alpha ,-11)^T}{\alpha \leq -7}$.
\subsection{Manufacturing system interpretation}
\label{subsec-machine}
\begin{figure}
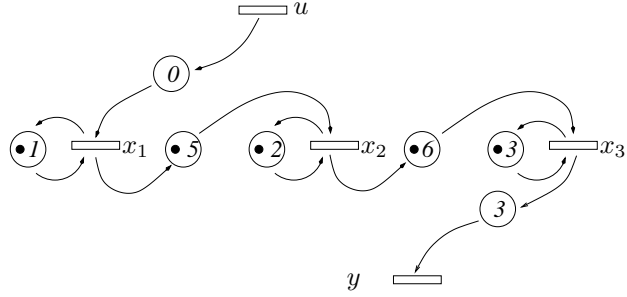

\begin{center}
\input fig7
\end{center}
\caption{A timed event graph representing three machines in tandem}
\label{fig-deds}
\end{figure}

We next interpret the previous computations
in terms of discrete event systems.
The dynamical system~\eqref{mpsystem},\eqref{exab} 
can be seen as the dater representation
of the 
timed event graph of Figure~\ref{fig-deds}
(we refer the reader to~\cite{bcoq} for more details on the
modeling of timed event graphs).

This graph represents three machines in tandem,
with respective processing times $1,2,3$. The first
machine is fed by a source $u$, and sends its
output to a second machine, with a transportation
delay of $5$. The second machine
sends its output to a third machine, with a transportation
delay of $6$. We associate
to each transition a dater function
$\N\to \rmax$: for instance, $u(k)$ gives
the date of the $k$-th firing of the transition
labeled $u$, i.e. $u(k)$ is the arrival
time of the $k$-th part, 
$x_3(k)$ gives the date at which the third
machine initiates its $k$-th task, etc.
The output transition labeled $y$ represents the times at which
finished parts become available.
The tokens in the places
$x_1\to x_2$ and $x_2\to x_3$ represent 
unfinished parts, which are initially available
when the activity starts.
The reachable spaces $\cR_k$, which were depicted
in Figure~\ref{figorth}, determine
the possible values of the daters $x_i(k)$, $i=1,\ldots,3$.
The assumption that $x(0)=\zero$ 
means that the machines are ready to operate
much before the first part arrives
from the source, so that only the $u\mapsto x$
relation is considered. 

Practically relevant quantities are the differences
$x_i-x_j$: for instance, $x_2(k)-x_1(k-1)$ gives
the sojourn time of the $(k-1)$-th part in the
storage resource between the first and the 
second machine. The timed event graph of 
Figure~\ref{fig-deds} is an example of instable system,
since the second machine is slower than
the first machine, parts may accumulate
infinitely in the intermediate storage resource.
This is reflected by the unboundedness of
the orthogonal representation of the reachable semimodule,
in Figure~\ref{figorth}. For instance, 
one can obtain for the trajectory $x(k)$ the
sequence of columns of the matrix $\cR_{\omega}$
by taking $u(k)=k-1$. 

The finite control sequence which leads
$x(0)=\zero$ to $x(k)=z$ can be computed
by solving the system $\Rmat_k U_k =z$,
where $U_k=(u(k),\ldots,u(1))^T$.
This system can be solved in polynomial time
using residuation theory, see~\cite{bcoq}
(or~\cite{butkovip94} for a more combinatorial
presentation).
We did not address the difficulty that 
the finite control sequence $U_k$ which leads to a given point of the
reachable space need not be physically admissible,
because dater functions must be {\em nondecreasing}.
Nondecreasing controls
can be modeled at the price of adding one variable:
if $u$ is an arbitrary control sequence, 
the max-plus linear
dynamical system $v(k)=v(k-1)\oplus u(k)$
computes the nondecreasing hull $v$
of $u$, and therefore $v$ represents
an arbitrary nondecreasing control sequence.

Let us now interpret the observable congruence $\cO_{\omega}$ in
terms of discrete event systems.
Consider two trajectories
$\{(x(k),y(k))\}_{k\geq 0}$ and
$\{(x'(k),y'(k))\}_{k\geq 0}$
associated to the same input sequence $\{u(k)\}_{k\geq 1}$. We can write
\begin{align}
y(k)=CA^kx(0)\oplus CA^{k-1}Bu(1)\oplus \cdots \oplus CBu(k)\enspace.
\label{output}
\end{align}
Comparing~\eqref{output} with the similar formula for $y'(k)$ we get that the following three assertions are equivalent:
\begin{enumerate}
\item\label{i1} the outputs $y$ and $y'$ corresponding to the zero input sequence satisfy 
$y(m)=y'(m)$ for all $0\leq m <k$,
\item\label{i2} for all input sequences, the associated outputs  satisfy 
$y(m)=y'(m)$ for all $0\leq m <k$,
\item\label{i3} $(x(0),x'(0))\in \cO_k$.
\end{enumerate}
In a ring, since addition is cancellative, the above assertions are equivalent to the following one:
\begin{enumerate}
\setcounter{enumi}{3}
\item the associated outputs $y$ and $y'$ corresponding to some input sequence $u$ satisfy 
$y(m)=y'(m)$ for all $0\leq m <k$.
\label{i4}
\end{enumerate}
The implication \ref{i4}$\Rightarrow $\ref{i2} is no longer true for linear systems over $\rmax$. In the sequel, we shall say that two initial conditions $x(0)$ and $x'(0)$ 
cannot be distinguished by observation up to time $k\in \N \cup \{\omega \}$ if any of the properties \ref{i1}--\ref{i3} holds. When $k=\omega$, we will simply say that $x(0)$ and $x'(0)$ cannot be distinguished by observation.

The congruence
$\cO_{\omega}$ obtained for the
transposed dynamics $A^T$ and observation matrix
$B^T$, already depicted in Figure~\ref{fig-dual}, 
corresponds to a timed
event graph in which the arcs are {\em reversed},
by comparison with the timed event graph of Figure~\ref{fig-deds}.

\begin{figure}
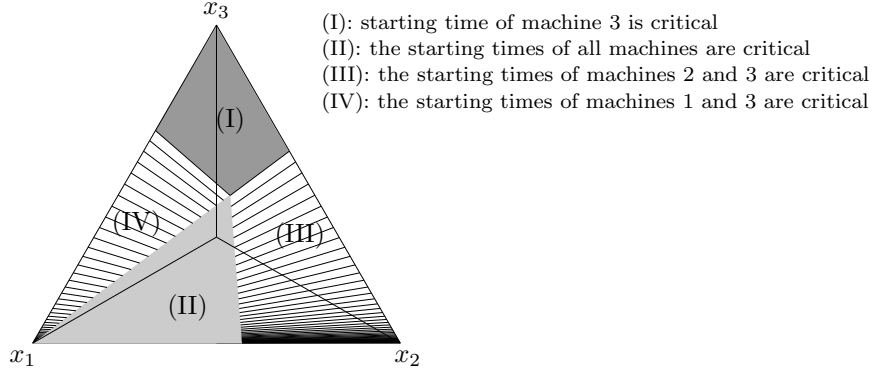

\begin{center}
\input fig8
\end{center}
\caption{Exponential projection of the observable congruence associated to the timed event graph of Figure~\ref{fig-deds}.}
\label{fig-obs}
\end{figure}

To give another example, with a more interesting
physical interpretation, let us introduce the 
observation matrix $C=(-\infty ,-\infty ,3)$,
which corresponds to the output $y=3x_3$ visible on Figure~\ref{fig-deds},
and consider the observable congruence $\cO_{\omega}$ corresponding
to the pair $(A,C)$, namely, $\cO_{\omega}=\ker \Omat_{\omega}$,
where
\[
\Omat_\omega=
\begin{pmatrix}
-\infty& -\infty&    3 \\
-\infty&   9 &    6 \\
14&  12&   9 \\
17& 15&  12\\
20& 18&  15\\
& \vdots
\end{pmatrix}
 \enspace .
\]
We have depicted in Figure~\ref{fig-obs} the
observable
congruence $\cO_\omega$ associated to $(A,C)$,
which is not only rational, but also finitely
generated (as a semimodule).
The figure represents four different types of equivalence classes associated
to finite points: for example the equivalence class 
of any point of the form $(x_3 -5,x_3
-3,x_3)^T$ is $\set{(x_1,x_2,x_3)^T}{x_1\leq x_3 -5,x_2\leq x_3 -3}$,
which is represented by the darkest tetragonal region,
labeled (I), in Figure~\ref{fig-obs}.
The light gray triangle in Figure~\ref{fig-obs}, which is labeled
(II), represents the
set of finite points of the semimodule over the min-plus semiring
$\zmin$ generated by the opposite of the rows of $\Omat_\omega$: we use again
the argument of~\cite{ccggq99} showing that equivalence classes for
$\cO_\omega$ of points in the interior of this set 
are singletons. Other equivalence classes are either half-lines,
as shown on the figure, or the singleton $\{(-\infty, -\infty,-\infty)\}$,
which cannot be represented by this projection.
Observation means looking at output times
of parts.
 Writing~\eqref{output} explicitly:
\begin{align}
\begin{array}{lcl}
y(0)&=&3x_3(0)\\
y(1)&=&9x_2(0)\oplus 6x_3(0)\\
y(2)&=&14x_1(0)\oplus 12 x_2(0)\oplus 9x_3(0) 
\end{array}
\end{align} 
and since the rows $4$, $5$, $6$, $\ldots $ of the observability matrix are proportional 
to the third row, it follows that two initial conditions cannot be distinguished by observation if and only 
if they lead to the same values for $y(0)$, $y(1)$ and
$y(2)$
(note that the input sequence 
can only change the values of $x_3(k)$ for $k\geq 3$). If we fix some
values for $x_3(0)$, 
$x_3(1)$ and $x_3(2)$, then we can determine the set of the initial 
conditions $(x_1(0),x_2(0),x_3(0))^T$ which will lead to these values and
therefore cannot be distinguished by observation. 
For example the set of
initial conditions which lead to the values $x_3(0)$, $x_3(1)=x_3(0)+3$
and $x_3(2)=x_3(0)+6$ is 
$\{(x_1(0),x_2(0),x_3(0))^T:x_1(0)\leq x_3(0)-5,x_2(0)\leq x_3(0) -3\}$,
that is, the equivalence class for
$\cO_\omega$ of the point $(x_3(0)-5,x_3(0)-3,x_3(0))^T$
(region (I)).
Therefore, this equivalence class represents a set of initial
conditions for which the starting time of machine $3$ is {\em critical}, 
which means that the output times are determined by this starting time. 
Similarly, the light grey region (II)
corresponds to a set of initial conditions 
such that $x_3(0)\leq 3+x_2(0)$ and $x_2(0)\leq 2+x_1(0)$. 
Every $x(0)$ in the interior of this zone, whose
equivalence class is a singleton, 
is such that the starting times of all machines are critical:
the output time of the first
part will only depend of the starting time
of machine $3$, the output time of the second part
will only depend of the starting time of machine
$2$, and the output time of the third part
will only depend of the starting time of machine
$1$. The half lines (III) and (IV) have a similar interpretation,
as summarized on Figure~\ref{fig-obs}.

\subsection{Rational semimodules over $\protect\rmax$ need not be stable by direct image}
\label{subsec-cex}
When $\cS=\rmax$,
the set of rational semimodules has no nice closure properties.
We first show that the direct image of a rational set
by a linear map need not be a rational set.
Let $\alpha$ denote a positive irrational
number, and consider the two vectors
\[
u= \begin{pmatrix}1\\-\alpha 
\end{pmatrix},
\qquad
v= -\alpha^{-1}u =\begin{pmatrix}-\alpha^{-1}\\ 1
\end{pmatrix}
\enspace,
\]
together with $R=\{u,v\}^*\setminus \{0\}$.
Since $R= \{u,v\}+\{u,v\}^*$,
$R$ is rational.  
Now, let $A=(0,0) \in (\rmax)^{1\times 2}$.
We have that:
\begin{align}
A(R) &=\set{\max (h_1-\alpha^{-1}h_2,h_2-\alpha h_1)}{h_1,h_2 \in \N, h_1+h_2\geq 1}\nonumber\\
&= \set{\max(-\alpha^{-1}t, t)}{t = h_2- \alpha h_1,\; h_1,h_2\in \N, 
h_1+h_2\geq 1} \enspace.\label{myAR}
\end{align}
We claim that the set $A(R)$ is not rational.
Indeed, let us assume by contradiction
that $A(R)$ is semilinear, i.e.,
that $A(R)=\cup_{1\leq i\leq k} (\{a_i\}+U_i^*) $, 
where the $a_i$ are elements of $\Rm$ and the $U_i$ are finite subsets
of $\Rm$. We first remark that since $-\infty\not\in A(R)$,
$a_i\neq-\infty$, and $-\infty\not\in U_i$,
for all $1\leq i\leq k$. Using this remark, we now deduce that
the elements of $U_i$ must be nonnegative:
otherwise, $A(R)$ would not be bounded from
below, and this would contradict the fact that $\inf A(R)=0$
which follows from~\eqref{myAR}.
Since all the elements of $U_i$ are nonnegative,
$A(R)$ has a minimal element (namely $\min_{1\leq i\leq k} a_i$),
and this contradicts~\eqref{myAR} because $\alpha$ is an
irrational number.

We next show that when $\cS=\rmax$, the image
of a rational semimodule by a linear map need
not be a rational semimodule.
Consider
\[
u= \begin{pmatrix}1\\-\alpha \\ 0 
\end{pmatrix},
\qquad
v= -\alpha^{-1}u =\begin{pmatrix}-\alpha^{-1}\\ 1\\ 0 
\end{pmatrix}
\enspace,
\]
\[
R= \{u,v\}^*\setminus \{0\}
,\qquad 
A=\begin{pmatrix}0 & 0 & -\infty \cr
 -\infty & -\infty & 0\end{pmatrix}\enspace,
\]
and $\cX=\vect R$. 
Then $A(\cX)= \vect A(R)$ is spanned
by the vectors
\[ 
\begin{pmatrix} \max (h_1-\alpha^{-1}h_2,h_2-\alpha h_1) \cr 0 \end{pmatrix},
\qquad \mrm{for}\qquad h_1,h_2\in \N, h_1+h_2\geq 1 \enspace.
\]
To make $A(\cX)$ more explicit, let us observe
that for all real numbers $\gamma,\delta$,
\begin{align}
\smallvect\begin{pmatrix}\gamma & \delta\\ 0 & 0
\end{pmatrix}
&= \Bigset{\begin{pmatrix} x_1\\ x_2\end{pmatrix}\in \R^2}{ \min(\gamma,\delta) + x_2 \leq x_1  \leq x_2 + \max(\gamma,\delta)}\cup
 \{\begin{pmatrix} -\infty\\ -\infty\end{pmatrix}\}.
\label{my2dsemi}
\end{align}
It follows from~\eqref{my2dsemi}
that
\begin{align}
A(\cX)= \Bigset{\begin{pmatrix} x_1\\ x_2\end{pmatrix}\in \R^2}{x_1>x_2}
\cup \{\begin{pmatrix} -\infty\\ -\infty\end{pmatrix}\}
\enspace.
\label{e-newsemiA}
\end{align}
Now, a straightforward variant of the proof
of the irrationality of $A(R)$ that we gave above
shows that $A(\cX)$ is not a rational semimodule,
for, if $A(\cX)$ was spanned by a semilinear
set, the quantity $x_1 -x_2$ would attain its
infimum when $x\in A(\cX)\cap \R^2$,
whereas~\eqref{e-newsemiA} shows that this infimum, which
is equal to $0$, is not attained. 

Thus, when $\cS=\rmax$, the direct image
of a rational semimodule by a linear map need not be rational.
\subsection{Noncommutative reachable spaces need not be rational}
\label{subsec-nc}
Let us consider now a time varying version
of the max-plus linear system (\ref{mpsystem}),
in which~\eqref{mata} is replaced by
\begin{align}
x(k)=A(k)x(k-1)\oplus Bu(k)  \enspace,
\end{align}
where the matrix $A(k)$ can take any value in a finite set
$\{ A_1,\ldots ,A_r\} \subset (\zmax)^{n\times n}$. 

In order to characterize the reachable
space and to show that it need not be rational,
it is useful to introduce some classical
automata theoretical notation (see~\cite{berstelreut}).
Let $\Sigma= \{ a_1,\ldots ,a_r\}$ denote an 
alphabet of $r$ letters. Recall that 
the free monoid $\Sigma^*$ is the set
of finite words on $\Sigma$, equipped
with concatenation product. 
Let $\mu: \Sigma^*\to  (\zmax)^{n\times n}$
denote the unique morphism of monoids
which sends $a_i$ to $A_i$.
The \NEW{reachable space} $\cR$, that is, the set of
all possible values of $x(k)\in (\zmax)^n$,
the control sequence $u$ and the time $k$
being chosen arbitrarily, starting from $x(0)=\zero$,
is given by: 
\[
\cR=\vect (\mu(\Sigma^*)B)
\enspace ,
\]
where we represent by $\vect (\mu(\Sigma^*)B)$ the subsemimodule which is generated by the columns of the 
matrices $\mu(w)B$, for $w\in \Sigma^*$.

We next show that $\mu(\Sigma^*)$,
and a fortiori $\set{\mu(w)B}{w\in \Sigma^*}$,
need not be rational subsets of $(\zmax)^{n\times n}$
and $(\zmax)^{n\times p}$, respectively,
and that the reachable space $\cR$ need
not be rational, a result which
illustrates a general difficulty of max-plus linear
semigroups (in a further work~\cite{gk02}, we show
that we cannot decide whether a matrix belongs to $\mu(\Sigma^*)$,
or whether a vector belongs to $\set{\mu(w)B}{w\in \Sigma^*}$).
In this paper, we will give a simple
counter-example, which relies on a remarkable construction of
I.~Simon~\cite{simon90}. To minimize
changes by comparison to~\cite{simon90},
we will work in the semiring $\zmin$,
rather than in $\zmax$. All the
results that follow have of course
equivalent versions in $\zmax$.

Let  $\nu: \{a_1,a_2\}^*\to (\zmin)^{4\times 4}$
denote the unique morphism such that:
\[\nu(a_1) =
\begin{pmatrix}
0 & \infty & \infty & \infty \cr 
\infty & 1 & 1 & \infty \cr 
\infty & \infty & \infty & \infty \cr
 \infty & \infty & \infty & 0 \end{pmatrix} 
\quad
\mrm{and} \quad
\nu(a_2)=\begin{pmatrix}
1 & 1 & \infty & \infty \cr
 \infty & \infty & \infty & 0 \cr
 \infty & \infty & \infty & 0 \cr
 \infty & \infty & \infty & 0 \end{pmatrix} 
\enspace ,\]
and consider the function $s: \{a_1,a_2\}^*\to \zmin, w\mapsto s(w)$,
\begin{align}
\label{ser-sim}
s(w)=\alpha\mu(w)\beta\quad\mrm{where}\;
\alpha=\begin{pmatrix} 0 & \infty & \infty & \infty \end{pmatrix}
\mrm{ and }
\beta=
\begin{pmatrix} 0 &\infty & \infty & 0 \end{pmatrix}^T\enspace. 
\end{align}
Simon~\cite{simon90} shows that
\begin{align}
\min \set{|w|}{s(w)\geq n} 
=\frac{n^2+n}{2} 
,\hspace{1em} \forall n\in \N\enspace ,
\label{growth}
\end{align}
where $|w|$ denotes the length of the word $w$.
In essence, \eqref{growth} means
that $s(w)$ takes values of order $\sqrt{|w|}$
when $|w|\to\infty$. We will use this property
to build an irrational reachable space $\cR$.

Let 
\[
D=\begin{pmatrix}
-1 & \infty\\
\infty & 0
\end{pmatrix}
\]
and consider the unique morphism
$\mu: \{a_1,a_2\}^*\to (\zmin)^{6\times 6}$,
\begin{align}
\label{def-mu}
\mu(a_1)=\diag(\nu(a_1),D),
\;\;
\mu(a_2)=\diag(\nu(a_2),D)
\enspace,
\end{align}
where $\diag(F,G)$ denotes the
matrix with diagonal blocks
$F$ and $G$ and $\infty$ elsewhere.
The following proposition shows that
the reachable space $\cR$ obtained
by taking
\begin{align}
\label{ser-sim2} 
B=\begin{pmatrix} 0 &\infty & \infty & 0 & 0 & 0 \end{pmatrix}^T\enspace. 
\end{align}
and $\mu$ as above, is irrational.

\begin{proposition}\label{prop-sp}
Let $\mu$ be defined by~\eqref{def-mu} and $B$ by~\eqref{ser-sim2}.
Then, the reachable space $\cR=\vect (\mu(\{a_1,a_2\}^*)B)$
is an irrational subsemimodule of $(\zmin)^6$.
Moreover , the semigroup $\mu(\{a_1,a_2\}^*)$ is an irrational
subset of $((\Zp)^{6\times 6},+)$.
\end{proposition}
\begin{proof}
Let $C$ denote the map
$(\zmin)^6\to (\zmin)^3$,
which is defined by the matrix: 
\[
C=\begin{pmatrix}
0 & \infty & \infty & \infty & \infty & \infty \cr 
\infty & \infty &  \infty &  \infty & 0 &  \infty \cr 
\infty & \infty & \infty & \infty &  \infty & 0 \end{pmatrix}. 
\]
Then we get that
\begin{align}
C\mu(w)B= (s(w),-|w|, 0)^T,\qquad \forall w\in \{a_1,a_2\}^*
\enspace .
\end{align}
If $\cX=\vect(\mu(\{a_1,a_2\}^*)B)$ were rational, 
$C(\cX)$ would also be rational, by Theorem~\ref{closureth}. 
We have represented $C(\cX)$ on Figure~\ref{fig-cs}:
the irrationality of $C(\cX)$ is intuitively clear
from the figure, since the boundary of the semimodule
has a discrete quadratic shape (extremal points
are represented by bold points).
\begin{figure}
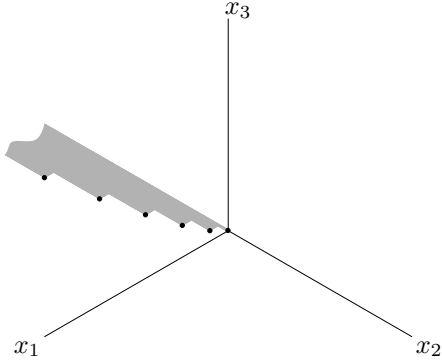

\begin{center}
\input fig9
\end{center}
\caption{An irrational subsemimodule of $(\zmin)^3$.}
\label{fig-cs}
\end{figure}

However, proving that the figure is correct
would require some reworking of the arguments
of~\cite{simon90}, so we will give a simpler
formal argument showing the irrationality
of $C(\cX)$. Since several words $w$ of the same
length can give the same $s(w)$,
the family $\{C\mu(w)B\}_{w\in \{a_1,a_2\}^*}=
\{(s(w),-|w|,0)^T\}_{w\in \{a_1,a_2\}^*}$
contains repeated elements.
So let us introduce a subfamily,
$\{C\mu(w)B\}_{w\in W}$,
with the property that 
$\set{C\mu(w)B}{w\in\{a_1,a_2\}^*}=
\set{C\mu(w)B}{w\in W}$, 
and
\begin{align}
\label{injec}
w,z\in W,\;|w|=|z|\implies s(w)\neq s(z)
\enspace. 
\end{align}
Let 
\[
W'= \set{w\in W}{(z\in W\setminus \{w\}\mrm{ and } s(z)\geq s(w))\implies |z|>|w|}
\enspace .
\]
We claim that 
\begin{align}
\forall w\in W',\;
C\mu(w)B \; \mrm{is an extremal point of } C(\cX).
\label{claim}
\end{align}
To show this, it suffices to check that
there is no family
$\{\lambda_z\}_{z\in W\setminus w}\subset\zmin$
such that
\[
C\mu(w)B=\bigoplus_{z\in W\setminus w} \lambda_z \otimes C\mu(z)B
\enspace ,
\]
i.e.
\begin{align}
\label{res}
(s(w),-|w|,0)^T=\inf_{z\in W\setminus w} \lambda_z+ (s(z),-|z|,0)^T
\enspace .
\end{align}
It follows from~\eqref{res} that
\[
\lambda_z \geq \max(s(w)-s(z),|z|-|w|,0) \enspace .
\]
Now, by definition of $W'$,
$\max(s(w)-s(z),|z|-|w|)>0$
for all $z\in W$ such that $z\neq w$, and since 
$\lambda_z\geq\max(s(w)-s(z),|z|-|w|)>0$
is an integer, we conclude that $\lambda_z\geq 1$.
Since this holds for all $z\in W\setminus w$,
the equality~\eqref{res} cannot hold,
because the third coordinate of the right-hand side of~\eqref{res}
must be greater than or equal to $1$, whereas
the third coordinate of the left hand side of~\eqref{res}
is equal to $0$. This shows~\eqref{claim}.

We finally show that $C(\cX)$ is irrational.
Consider
\begin{align}
E=\set{(x_1,x_2)}{(x_1,x_2,0)^T\in \Ext(C(\cX))}
\enspace .
\label{extremal}
\end{align}
If $\cX$ were rational, so would be $C(\cX)$,
and by Theorem~\ref{extrat},
the set of extremal points $\Ext(C(\cX))$
of $C(\cX)$ would be rational, and so $E$
would be rational.

Now, it follows from the definition
of extremal points that for any set 
$G$ of generators of a semimodule $\cX$,
\begin{equation}
\Ext(\cX)\subset \Z+G=\set{\lambda\otimes g}{\lambda\in \Z,\;g\in G}
\enspace.
\label{wherext}
\end{equation}
Combining~\eqref{wherext},~\eqref{extremal}, and~\eqref{claim},
and using the fact that the third coordinate
of $C\mu(w)B$ is $0$ for all $w\in \{a_1,a_2\}^*$,
we get that 
\begin{align}
\label{sandwitch}
\set{(s(w),-|w|)}{w\in W'} \subset
E \subset
\set{(s(w),-|w|)}{w\in \{a_1,a_2\}^*}
\enspace .
\end{align}
Now, for any rational subset $R$ of $(\Z^2,+)$,
consider the function:
\[
\gamma_R: \Z\to \Z\cup\{\pm\infty\},\;
\gamma_R(n)=\sup\set{k\in \Z}{(n,k)\in R}
\enspace ,
\]
together with its support:
\[
\mrm{supp}\,\gamma_R=\set{n\in\Z}{\exists k\in \Z,\;(n,k)\in R}
=\set{n\in \Z}{\gamma_R(n)\neq-\infty} \enspace .
\]
It follows from the fact
that rational subsets of $(\Z^2,+)$ are semilinear
that if $R$ is rational, the restriction of $\gamma_R$
to its support can be bounded
from below by an affine function
when $n\to\infty$.
But~\eqref{sandwitch} together with~\eqref{growth}
show that $\gamma_E(n)=-(n^2+n)/2$. Therefore, $E$ is
irrational, a contradiction.
\end{proof}
The counter example of Proposition~\ref{prop-sp}
shows that the rational semimodules tools
do not apply naturally to max-plus automata
problems, such as the ones
appearing in~\cite{gaubert95c,klimann,klimann99}.

\medskip\noindent{\em Acknowledgement.}\/
The authors thank Guy Cohen and Jean-Pierre Quadrat
for having inspired this work by many useful discussions.

\end{document}